\documentclass[sn-mathphys-num]{sn-jnl}

\usepackage{graphicx}%
\usepackage{multirow}%
\usepackage{amsmath,amssymb,amsfonts}%
\usepackage{amsthm}%
\usepackage{mathrsfs}%
\usepackage[title]{appendix}%
\usepackage{xcolor}%
\usepackage{textcomp}%
\usepackage{manyfoot}%
\usepackage{booktabs}%
\usepackage{algorithm}%
\usepackage{algorithmicx}%
\usepackage{algpseudocode}%
\usepackage{listings}%
\usepackage{hyperref}

\usepackage{tikz, pgfplots}
\pgfplotsset{compat=1.15}

\tikzstyle{mystyle}=[line join=round, line cap = round, line width = 0.5pt, >=latex, font=\footnotesize, x=0.85cm,y=0.85cm]
\tikzstyle{mylatstyle}=[mystyle,x=0.55cm,y=0.55cm]
\usetikzlibrary{arrows}

\newcommand{\lattice}[2]{\foreach \i in {0,1,...,#1}{
		\foreach \j in {0,1,...,#2}{
			\draw[fill] (\i,\j) circle (0.025cm);}}}

\def\do#1{\csdef{#1}{\mathbb{#1}}}
\docsvlist{N,Z,Q,R,C}

\newcommand{\conv}{\operatorname{conv}}

\theoremstyle{thmstyleone}%
\newtheorem{theorem}{Theorem}
\newtheorem*{theorem*}{Theorem}
\newtheorem{lemma}[theorem]{Lemma}

\theoremstyle{thmstyletwo}%
\newtheorem{example}{Example}%

\theoremstyle{thmstylethree}%

\begin{document}
\title[Article Title]{Convex Lattice Polygons with $k\ge3$ Interior Points}


\author*[1]{\fnm{Dana} \sur{Paquin}}\email{dpaquin@calpoly.edu}

\author[1]{\fnm{Elli} \sur{Sumera}}\email{esumera@calpoly.edu}

\author[1]{\fnm{Tri} \sur{Tran}}\email{htran82@calpoly.edu}


\affil*[1]{\orgdiv{Department of Mathematics}, \orgname{California Polytechnic State University}, \orgaddress{\street{1 Grand Avenue}, \city{San Luis Obispo}, \postcode{93407}, \state{California}, \country{USA}}}






\abstract{We study the geometry of convex lattice $n$-gons with $n$ boundary lattice points and $k\geq 3$ collinear interior lattice points. We describe a process to construct a primitive lattice triangle from an edge of a convex lattice $n$-gon, hence adding one edge in a way so that the number of boundary points increases by $1$, while the number of interior points remains unchanged. We also present the necessary conditions to construct such a primitive lattice triangle, as well as an upper bound for the number of times this is possible. Finally, we apply the previous results to fully classify the positive integers for which there exists a convex $n$-gon with $k$ collinear and non-collinear interior points.}





\maketitle

\section{Introduction}\label{sec:1}

A \textit{lattice point} $(x,y)$ is an element of $\Z\times\Z = \Z^2$. A \textit{lattice polygon} is a simple polygon whose vertices are lattice points. By the Jordan Curve Theorem (cf. \cite{jordan,tverberg}), $P$ has a well-defined interior and exterior, denoted $\operatorname{int}(P)$ and $\operatorname{ext}(P)$, respectively. If $P$ is a lattice polygon, we let $B(P)$ and $I(P)$ denote the number of boundary and interior lattice points of $P$, respectively. If the context is clear, we may omit the parentheses and simply write $B$ or $I$. Additionally, let $L(v,w)$ denote the number of lattice points on the line segment connecting two points $v,w\in\Z^2$. We can compute the area of any lattice polygons simply by counting its boundary and interior points, that is,
\begin{align*}
	\operatorname{area}(P) = \frac{B(P)}{2} + I(P) - 1.
\end{align*}
This is the well-known Pick's theorem (cf. \cite{pick}). 

Two lattice points $v,w\in\Z^2$ are \textit{visible} if $L(v,w)=0$. It was shown in \cite[Theorem 3.8]{apostol} that two points $v = (x_1,y_1)$ and $w = (x_2,y_2)$ are visible if and only if $\gcd(x_1-x_2,y_1-y_2) = 1$. Unless otherwise stated, we assume all pairs of adjacent lattice points on any polygon are visible. Equivalently, $B(P) = n$ for all $n$-gons $P$.

There is a lot of interest regarding the geometric structure of lattice polygons having a certain number of interior points (cf. \cite{coleman,kolodziejczyk1,kolodziejczyk2,li,rabinowitz1,rabinowitz2,rabinowitz3,simpson,wei}). We restrict our focus on lattice polygons with collinear interior points. A \textit{$2$-collinear} integer $k$ is a positive integer such that some lattice triangle with $B=3$ has $k$ interior points, then the interior points are necessarily collinear. Li and Paquin showed in \cite{li} that $1$, $2$, $4$, and $7$ are the only $2$-collinear integers. If we omit the values of $k$ where the interior points are trivially collinear, i.e., $1$ and $2$, then $2$-collinearity is a property intrinsic to lattice triangles. In other words, for $n\ge4$ there does not exist an interger $k>2$ in which if an $n$-gon with $B = n$ has $I = k$, then the interior points must be collinear.

A lattice polygon is \textit{convex} if for every pair of adjacent vertices $v_i,v_j$, the line
\begin{align*}
	(1-t)v_i+tv_j,\quad t\in(0,1)
\end{align*}
is fully contained in its interior. When appropriate we shall use the alternate definition from \cite[Lemma 1.18]{orourke} that a lattice polygon is convex if and only if all interior angles are strictly less than $\pi$. For any list $v_1,\ldots,v_n\in\Z^2$, let $\conv\{v_1,\ldots,v_n\}$ denote its \textit{convex hull}. For brevity we assume all polygons in this paper are convex, unless otherwise specified. We state the main results of the paper.

\begin{theorem*}
	The only integers $n$ in which we can construct a convex $n$-gon with $n$ boundary and $k\ge3$ interior points in which the interior points are collinear are $3$, $4$, $5$, and $6$.
\end{theorem*}

\noindent This result means that there is a maximum number of edges for an $n$-gon with collinear interior points, but there exists such a polygon for all $k\ge3$ interior points. Conversely, the next result states that we can also construct a polygon with the same number of edges but the interior points need not be collinear. 

\begin{theorem*}
	For any $k\ge4$ there exists a convex $n$-gon, $n\in\{4,5,6\}$ such that the $k$ interior points are not collinear.
\end{theorem*}



\subsection{An Equivalence Relation on Lattice Polygons}

An \textit{affine transformation} is a map $F:\R^2\to\R^2$ defined by $v\mapsto Av+b$ for $A\in\operatorname{GL}_2(\R)$ and $b\in\R^2$. It is \textit{integral} if $A\in\operatorname{GL}_2(\Z)$ and it is \textit{unimodular} if $|\det A| = 1$, or preserves area. It was shown in \cite[Theorem 15.1]{martin} that affine maps preserve collinearity. We state a more general result than \cite[Lemma 2.3]{li}.

\begin{lemma}
	Suppose $F:\R^2\to\R^2$ is an integral unimodular affine transformation and $P$ a lattice polygon. Then $B(F(P)) = B(P)$, $I(F(P)) = I(P)$, and $F(P)$ is convex if and only if $P$ is convex. 
\end{lemma}


\begin{proof}
	Denote the vertices of $P$ by $v_1,\ldots,v_n$ in counter-clockwise order. If $\gamma:[0,1]\to \R^2$ is a simple closed curve, then $F\circ\gamma$ is also simple because otherwise there exist $t_1,t_2\in (0,1)$, $t_1\neq t_2$ such that $F\circ\gamma(t_1) = F\circ\gamma(t_2)$, contradicting injectivity of $F$ since $\gamma(t_1)\neq\gamma(t_2)$. Moreover, $\gamma(0) = \gamma(1)$ since it is closed, so $F\circ\gamma(0) = F\circ\gamma(1)$, showing $F\circ\gamma$ is also closed. Since any lattice polygon is a simple closed curve, its image under $F$ must also be simple and closed, hence a lattice polygon.
	
	We compute
	\begin{align*}
			B(F(P)) &= n + L(F(v_n),F(v_1)) + \sum_{i=1}^{n-1}L(F(v_{i}),F(v_{i+1}))\\
			&= n + L(v_n,v_1) + \sum_{i=1}^{n-1}L(v_{i},v_{i+1})\\
			&= B(P),
		\end{align*}
	and since $F$ preserves area,
	\begin{align*}
			I(F(P)) &= \operatorname{area}(F(P)) - \frac{B(F(P))}{2} + 1\\ 
			&= \operatorname{area}(P) - \frac{B(P)}{2} + 1\\
			&= I(P).
		\end{align*}
	
	Finally, pick any non-adjacent vertices $v_i,v_j\in\partial P$. Then $(1-t)v_i + tv_j\in P$ for all $t\in[0,1]$. But then this is true if and only if $F((1-t)v_i + tv_j) = (1-t)Fv_i + tFv_j \in F(P)$ for all $t\in[0,1]$. Since $Fv_i,Fv_j\in\partial F(P)$, $(1-t)Fv_i + tFv_j\in\operatorname{int}(P)$ for all $t\in(0,1)$. Therefore if $P$ is convex then so is $F(P)$ and vice versa.
\end{proof}

From this result we can define an equivalence relation on all lattice polygons by asserting two polygons are related, or equivalent, if and only if there exists an integral unimodular affine transformation $F$ that sends one to the other. This greatly simplifies our construction of lattice polygons as we may pick a representative that is easy to work with. In Section \ref{sec:4.2} we show that for polygons with collinear interior points, it suffices to consider those with interior points lying on the $x$-axis.

\subsection{Isometries}

For some of the computations in Section \ref{sec:4}, there may not exist an integral unimodular affine map that satisfies our requirement as the angle required to rotate an edge of a polygon into the correct subset of a quadrant might not be a lattice angle. Instead we consider affine transformations that preserve distance and angles, but not necessarily integral, that is, it need not send lattice points to lattice points. 

An \textit{isometry} is a linear map that preserves distance between two vectors. Here we are working in the plane, so we use the standard Euclidean distance. It is known that isometries preserves Euclidean angle and has determinant $\pm1$. To be more specific, we refer to the following result.

\begin{theorem}[{\cite[Theorem 8.6]{martin}}]
	Every non-identity isometry of $\R^2$ is either a translation, rotation, reflection, or glide reflection.
\end{theorem}

\section{Appending Primitive Triangles}\label{sec:4}

A \textit{primitive lattice polygon} is a lattice $n$-gon with $B = n$ and $I = 0$. From Pick's theorem it can be shown that $T$ is a primitive lattice triangle if and only if $\operatorname{area}(T) = \frac{1}{2}$. 

%
%

\begin{figure}
	\centering
	\begin{tikzpicture}[mystyle]
		\draw[->] (-2,0) -- (5,0);
		\draw[->] (0,-1.5) -- (0,4);
		\draw[thick] (0,0) -- node[below right] {$L_i$} (3,2.5);
		\draw[thick,dashed] (0,0) -- (-2,-1) node[below] {$L_{i-1}$};
		\draw[thick,dashed] (3,2.5) -- (4,4) node[right] {$L_{i+1}$};
		
		\draw (0.75,0) arc (0:39.8:0.75) node[pos=0.5,right] {$\theta_i$};
		\draw (0.45,0) arc (0: 206.565:0.45) node[pos=0.7,above left] {$\theta_{i-1}$};
		\draw (3.5,2.5) arc (0:56.31:0.5) node[pos=0.65,right] {$\theta_{i+1}$};
		
		\draw[dotted] (0,0) -- (3,3) node[above left] {$y=x$};
		\draw[shift={(3,2.5)}] (0,0) -- (1,0);
		\node at (1,2.5) {$\operatorname{int}(P)$};
		\node at (3,1) {$\operatorname{ext}(P)$};
		\draw[fill] (0,0) circle (0.04cm);
		\draw[fill] (3,2.5) circle (0.04cm) node[below right] {$(p,q)$};
	\end{tikzpicture}
	\caption{The standard orientation of an edge $L_i$.} \label{fig:canon_rep}
\end{figure}
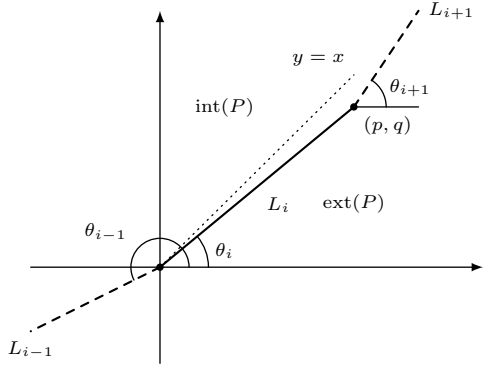

Given any polygon $P$, we fix an edge, translate the origin up to one of its vertices, and orient via rotation and reflection such that the interior of $P$ lies to the upper left and the exterior lies to the lower right of the edge. By assumption, there are no lattice points on the edge, so one vertex is $(0,0)$ and the other is $(p,q)$ with $\gcd(p,q) = 1$. Denote the line segment from $(0,0)$ to $(p,q)$ by $L_i$ and its adjacent edges $L_{i+1}$ and $L_{i-1}$, reading counter-clockwise. Let $\theta_{i}$, $\theta_{i+1}$, and $\theta_{i-1}$ denote the polar angle $L_i$, $L_{i+1}$, and $L_{i-1}$ make with the $x$-axis, respectively. Note that by this definition, $\theta_{i-1},\theta_i,\theta_{i+1}\in[0,2\pi]$. It is clear from geometry that $\theta_i = \arctan\frac{q}{p}$. Further, we may assume without loss of generality that $p>q$, so that the slope of the line is always less than $1$; otherwise we can always reflect and rotate the polygon so that the condition is satisfied. We call this the \textit{standard orientation}, as shown in Figure \ref{fig:canon_rep}.

It is important to reiterate that the affine transformation that sends the edge $L_i$ to its standard orientation must be unimodular, but not necessarily integral. This is because we are only concern about the slope of the edges and the fact that $(p,q)$ are coprime. Thus we only consider the image of $P$ under isometries, i.e., transformations that preserve lattice distances and angles. 

Before continuing with the technical results, we review a key concept from calculus. Let $L:ax+by+c=0$ be a line in $\R^2$. For any point $(x,y)\in\R^2$, the distance from $(x,y)$ to $L$ is given by \cite{spain}
\begin{align*}
	d_L(x,y) = \frac{|ax+by+c|}{\sqrt{a^2+b^2}}.
\end{align*}

\begin{lemma}\label{lem:existence_convex_m1}
	Let $P$ be a convex lattice $n$-gon with $B(P) = n$ and $I(P) = k$. Suppose $P$ has an edge $L_i$ with coordinates (up to isometry) $(0,0)$ and $(p,1)$, $p\in\N$. Further, assume
	$\frac{\pi}{2}<\theta_{i+1}<\theta_i+\pi$ if $p=1$ or $\arctan\frac{1}{p-1}<\theta_{i+1}<\theta_i + \pi$ if $p>1$, and $\theta_i<\theta_{i-1}<\pi$. Then there exists a convex lattice $(n+1)$-gon $P^{+}$ with $B(P^+) = n+1$ and $I(P^+) = k$.
\end{lemma}

\begin{proof}
	Note that the conditions $\theta_{i+1}<\theta_i+\pi$ and $\theta_{i-1}>\theta_{i}$ are implied from the supposition that $P$ is convex. Draw an edge from $(0,0)$ to $(1,0)$, from then draw another edge to $(a,b)=(p,1)$. Then $L((1,0),(a,b))=0$ since $\gcd(p-1,1)=1$. The triangle defined by these three vertices has area $\frac{1}{2}$, hence a primitive triangle. By removing the edge $L_{i}$, we have created an $(n+1)$-gon with $n+1$ boundary points and $k$ interior points. The condition $\theta_{i-1}<\pi$ guarantees that the vertices of $L_{i-1}$ and $(1,0)$ are not collinear. Finally, the condition $\theta_{i+1}>\frac{\pi}{2}$ (if $p=1$) or $\theta_{i+1}>\arctan\frac{1}{p-1}$ (if $p>1$) ensures the vertices of $L_{i+1}$ and $(1,0)$ are also not collinear, hence preserving convexity of $P$.
\end{proof}

\begin{lemma}\label{lem:existence_convex_mn}
	Let $P$ be a convex lattice $n$-gon with $B(P) = n$ and $I(P) = k$. Suppose $P$ has an edge $L_i$ with coordinates (up to isometry) $(0,0)$ and $(p,q)$, $p,q>1$, and $\gcd(p,q) = 1$. By B\'{e}zout's lemma there exists $a',b'\in\Z$ such that $|a'|$ is smallest and $a'q+b'p=\pm1$. Define 
	\begin{align*}
		(a,b) = \begin{cases}
			(a',-b') &\text{if }-\frac{b'}{a'}<\frac{q}{p},\\
			(p-a',q+b') &\text{if }-\frac{b'}{a'}>\frac{q}{p}.
		\end{cases}
	\end{align*}
	Further, assume $\arctan\frac{q-b}{p-a}<\theta_{i+1}<\theta_i+\pi$ and $\theta_i<\theta_{i-1}<\arctan\frac{b}{a} + \pi$. Then there exists a convex lattice $(n+1)$-gon $P^{+}$ with $B(P^+) = n+1$ and $I(P^+) = k$.
\end{lemma}

\begin{proof}
	Note that the conditions $\theta_{i+1}<\theta_i+\pi$ and $\theta_{i-1}>\theta_i$ are implied from the supposition that $P$ is convex. 
	
	Given the line segment $L_i = \overline{(0,0)(p,q)}$, with $p>q>1$ and $\gcd(p,q)=1$, we wish to find a unique point $(a,b)$ such that the distance between $(a,b)$ and $L_i$ is minimized. We parametrize $L_i:qx - py = 0$. For any point $(x,y)\in\Z^2$, the distance from $(x,y)$ to $L_i$ is 
	\begin{align*}
		d_{L_i}(x,y) = \frac{|xq - yp|}{\sqrt{p^2+q^2}}.
	\end{align*}
	By B\'{e}zout's lemma, there exists $a',b'\in\Z$ such that $a'>0$, $|a'|$ is smallest and $a'q+b'p=\pm1$. Here we relax the greatest common divisor to be $\pm1$ since we are measuring absolute distance, so sign does not matter. This forces $b'<0$ and $\frac{b'}{a'}<0$.
	
	If $-\frac{b'}{a'}<\frac{q}{p}$ then the point $(a',-b')$ lies outside $P$, so set $(a,b) = (a',-b')$. The distance is thus 
	\begin{align*}
		d_{L_i}(a,b) = \frac{|aq-bp|}{\sqrt{p^2+q^2}} = \frac{|a'q+b'p|}{\sqrt{p^2+q^2}} = \frac{1}{\sqrt{p^2+q^2}}.
	\end{align*}
	If $-\frac{b'}{a'}>\frac{q}{p}$ then $(a',-b')$ lies in the interior of $P$. But we can reflect it about $L_i$ by setting $(a,b) = (p-a',q+b')$. Note that $-\frac{b'}{a'}>\frac{q}{p}$ is equivalent to $-b'p>a'q$, or $b'p<-a'q$. Adding $pq$ on both sides and factor gives $p(q+b')<q(p-a')$, which is equivalent to $\frac{q+b'}{p-a'} = \frac{b}{a}<\frac{q}{p}$. Here $p-a'\neq0$ because otherwise $1=a'q+b'p=p(q+b')$, which implies $q+b' = \frac{1}{p}$, a contradiction since $p>1$ and $q+b'\in\Z$. (We can also prove $q + b'\neq0$ using the same argument.) 
	This shows $(a,b)$ lies outside $P$, and the distance 
	\begin{align*}
		d_{L_i}(a,b) &= \frac{|aq-bp|}{\sqrt{p^2+q^2}}\\
		&= \frac{|(p-a')q-(q+b')p|}{\sqrt{p^2+q^2}}\\
		&= \frac{|a'q-b'p|}{\sqrt{p^2+q^2}}\\
		&= \frac{1}{\sqrt{p^2+q^2}}
	\end{align*}
	is minimized. Lastly, if $-\frac{b'}{a'}=\frac{q}{p}$ then, since $\gcd(p,q)=1$, $\frac{q}{p}$ is in reduced form, so it follows $-b'=q$ and $a'=p$. But this means $\pm1=a'q+b'p=pq-pq=0$, a contradiction. Thus this case is not possible, and we omit from the conditions. In all cases, we have $a>0$ and $b>0$.
	
	It remains to address the possibility of $(a,b)$ lies outside $\conv\{(0,0),(p,0),(p,q)\}$. We covered the case $a = p$ in the previous paragraph. Assume $a>p$, then $a = p+j$ for some $j\in\N$. If $b<q$ then $0<q-b<q$. We get
	\begin{align*}
		1 &= |aq-bp|\\
		&= |(p+j)q - bp|\\
		&= |p(q-b) + jq|\\
		&> |p(q-b) + j(q-b)|\\
		&= (p+j)(q-b).
	\end{align*}
	Since $p+j>0$, it must be that $q-b<1$, or $q<1+b$, or $q\le b$, a contradiction. Thus $a>p$ forces $b\ge q$. Now $d_{L_i}(a,b) = \frac{1}{\sqrt{p^2+q^2}}$, but then 
	\begin{align*}
		d_{L_i}(a-p,b-q) = \frac{|(a-p)q-(b-q)p|}{\sqrt{p^2+q^2}} = \frac{|aq-pq-bp+pq|}{\sqrt{p^2+q^2}} = \frac{1}{\sqrt{p^2+q^2}}.
	\end{align*}
	Therefore $(a-p,b-q)$ is a lattice point inside the $\conv\{(0,0),(p,0),(p,q)\}$ with the same distance, which contradicts the minimality of $|a'|$. So we obtain $0<a<p$ and $0<b<q$. Lastly, we can see that $(a,b)$ is unique in $\conv\{(0,0),(p,0),(p,q)\}$ because otherwise there exists a point $(a_1,b_1)$, also in the triangle with $d_{L_i}(a_1,b_1) = d_{L_i}(a,b)$, so they both lie on the line parallel to $L_i$. But then this implies existence of a lattice point between $(0,0)$ and $(p,q)$, contradicting the relatively prime condition on $p$ and $q$.
	
	Consider the triangle $\conv\{(0,0),(a,b),(p,q)\}$. The base has length $\sqrt{p^2+q^2}$ and the height has length $d_{L_i}(a,b) = \frac{1}{\sqrt{p^2+q^2}}$. Thus the area is $\frac{1}{2}$, which implies the triangle is primitive. The condition $\theta_{i+1}>\arctan\frac{q-b}{p-a}$ guarantees that the vertices of $L_{i+1}$ and the point $(a,b)$ are not collinear, or breaking convexity. The condition $\theta_{i-1}<\arctan\frac{b}{a} + \pi$ guarantees that the vertices of $L_{i-1}$ and $(a,b)$ are also not collinear, or breaking convexity. Therefore appending this triangle to $P$ gives an $(n+1)$-gon with $n+1$ boundary and $k$ interior points.
\end{proof}

The final case that we have to check is $(p,q) = (1,0)$ or $(0,1)$. The only possibilities for primitive triangles to append are $(a,b) = (0,-1)$ or $(a,b) = (1,-1)$. Without loss of generality, consider the former case. Convexity gives $0<\theta_{i-1}<\frac{\pi}{2}$ and $\frac{\pi}{4}<\theta_{i+1}<\pi$. The only way for this to happen so that $P$ is a lattice polygon is that $L_{i-1}$ and $L_{i+1}$ intersect at $(1,1)$. But this forces $P$ to be a triangle. Therefore if $n\ge4$, we cannot append any primitive triangle to an edge of length $1$.

In fact, Lemmas \ref{lem:existence_convex_m1} and \ref{lem:existence_convex_mn} are biconditionals; that is, if $P$ is a lattice $n$-gon and it is possible to create (from $P$) an $(n+1)$-gon $P^{+}$ with $B(P^{+}) = n+1$ and $I(P^{+}) = I(P)$, then there exists an edge of $P$ that satisfies the hypotheses of one of the lemmas. This result does not contribute to proving the main result, so we omit the proof.

\begin{example}
	This example shows a polygon $P$ before and after having primitive triangles appended to the edges that satisfy the conditions of the technical lemmas. Note further that we cannot append a primitive triangle to any edge of the polygon on the right.
	
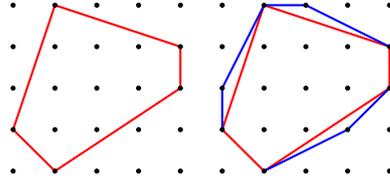
\begin{center}
		\begin{tikzpicture}[mylatstyle]
			\draw[thick,red,shift={(1,0)}] (0,0) -- (3,2) -- (3,3) -- (0,4) -- (-1,1) -- cycle;
			\begin{scope}[shift={(6,0)}]
				\draw[thick,red] (0,0) -- (3,2) -- (3,3) -- (0,4) -- (-1,1) -- cycle;
				\draw[thick,blue] (0,0) -- (2,1) -- (3,2);
				\draw[thick,blue] (3,3) -- (1,4) -- (0,4) -- (-1,2) -- (-1,1);
			\end{scope}
			\lattice{9}{4}
		\end{tikzpicture}
	\end{center}
\end{example}

A natural question arises: For a given lattice polygon $P$, how many times can we append a primitive triangle so that convexity is preserved and the number of interior points remains the same? More specifically, once we have attached a triangle to an edge, we have created two different edges; so can we append (possibly smaller) primitive triangles to the new edges? Intuitively we speculate that the process can only go so far before introducing a non-vertex boundary point or breaking convexity.

\begin{figure}
	\centering
	\begin{tikzpicture}[mystyle]
		\node (0) at (0,0) {};
		\node (1) at (5,4) {};
		\node (a) at (4,2) {};
		\node (l) at (2,0.75) {};
		\node (u) at (4.75,3) {};
		\draw[thick] (0.center) -- node[above left] {$L_i$} (1.center) -- node[left] {$u$} (a.center) -- node[pos=0.35,above] {$\ell$} cycle;
		\draw[blue,thick] (0.center) -- (l.center) -- (a.center);
		\draw[blue,thick] (a.center) -- (u.center) -- (1.center);
		
		\draw[fill] (0,0) circle (0.04cm) node[left] {$(0,0)$};
		\draw[fill] (1.center) circle (0.04cm) node[above] {$(p,q)$};
		\draw[fill] (a.center) circle (0.04cm) node[below right] {$(a,b)$};
		\draw[fill] (l.center) circle (0.04cm) node[below right] {$(a_{\ell},b_{\ell})$};
		\draw[fill] (u.center) circle (0.04cm) node[right] {$(a_u,b_u)$};
	\end{tikzpicture}
	\caption{Geometric interpretation of appending primitive triangles to the upper and lower edges of the previously appended primitive triangle.}
	\label{fig:2}
\end{figure}
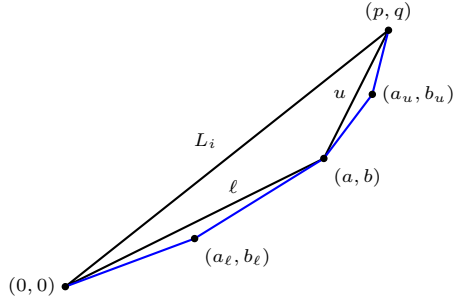

Suppose $P$ is a polygon with an edge such that we can apply Lemma \ref{lem:existence_convex_m1} or \ref{lem:existence_convex_mn}. Then the triangle $\conv\{(0,0),(a,b),(p,q)\}$ is primitive. So we get $\gcd(a,b)=\gcd(p-a,q-b)=1$. Applying B\'{e}zout's lemma again to obtain $(a_{\ell},b_{\ell})$ and $(a_u,b_u)$ as shown in Figure \ref{fig:2}. These two points are the closest lattice points to the edges $\ell$ and $u$, respectively.

\begin{theorem}\label{thm:at_most_n_triangles}
	For any convex $n$-gon $P$ with $B(P)=n$ and $I(P)=k$, we can append at most one primitive triangle to each edge without breaking convexity, or introducing non-vertex boundary points.
\end{theorem}

\begin{proof}
	It suffices to show that for any edge $L_i$ that satisfies the conditions for appending a primitive triangle, we can do so only once. Fix an edge $L_i$ of $P$. By a translation and isometric transformation, we may assume it has the standard orientation in Figure \ref{fig:canon_rep}. From here we compute the B\'{e}zout coefficients and apply the algorithm in Lemmas \ref{lem:existence_convex_m1} and \ref{lem:existence_convex_mn}.
	
	If $p=q=1$, then $(a,b) = (1,0)$, so both edges $u$ and $\ell$ have length $1$, so we cannot append a primitive triangle because the least amount of edges of a polygon after appending a triangle is $4$. If $q=1$ and $p>1$, then $(a,b) = (1,0)$, so we cannot append a primitive triangle to $\ell$. We can also see $(a_u,b_u) = (2,0)$, so $(0,0)$, $(a,b)$, and $(a_u,b_u)$ are collinear, also breaking convexity.
	
	Assume $(p,q)$ is such that $p>q>1$. We prove that it is not possible to append a primitive triangle to $\ell$. Suppose we have computed $(a,b)$. Draw a line $L_i'$ parallel to $L_i$ that passes through $(a,b)$, then draw a line $L_i''$ connecting $(a,b)$ and $(p,q)$. This is shown in the following figure. Here the red triangle is the primitive triangle $\conv\{(0,0),(a,b),(p,q)\}$.
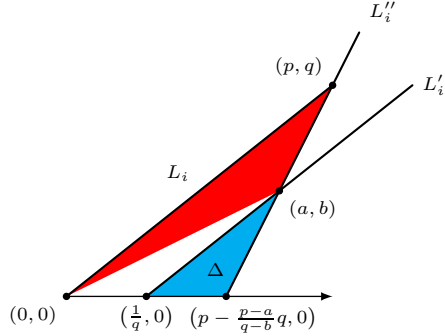
	\begin{center}
		\begin{tikzpicture}[mystyle,x=0.7cm,y=0.7cm]
			\fill[red,opacity=0.25] (0,0) -- (5,4) -- (4,2) -- cycle;
			\fill[cyan,opacity=0.25] (1.5,0) -- (4,2) -- (3,0) -- cycle;
			
			\draw[->] (0,0) -- (5,0);
			\draw[thick] (0,0) -- node[above left] {$L_i$} (5,4);
			\draw[thick] (1.5,0) -- (6.5,4) node[right] {$L_i'$};
			\draw[thick] (3,0) -- (5.5,5) node[above right] {$L_i''$};
			
			\draw[fill] (0,0) circle (0.04cm) node[below left] {$(0,0)$};
			\draw[fill] (5,4) circle (0.04cm) node[above left] {$(p,q)$};
			\draw[fill] (4,2) circle (0.04cm) node[below right] {$(a,b)$};
			\draw[fill] (1.5,0) circle (0.04cm) node[below] {$\big(\frac{1}{q},0\big)$};
			\draw[fill] (3,0) circle (0.04cm) node[below,shift={(0.5,0)}] {$\big(p-\frac{p-a}{q-b}q,0\big)$};
			\node at (2.8,0.5) {$\Delta$};
		\end{tikzpicture}
	\end{center}
	Previously we showed $d_{L_i}(x,y) = \frac{1}{\sqrt{p^2+q^2}}$ for all $(x,y)\in L_i'$. By construction, $(a_{\ell},b_{\ell})$ must lie outside the red triangle. It cannot lie on or to the right of $L_{i}''$ because it would break convexity. If $(a_{\ell},b_{\ell})$ lies in or on the boundary of $\conv\big\{(0,0),\big(\frac{1}{q},0\big),(a,b)\big\}$, then $d_{L_i}(a_{\ell},b_{\ell})\le d_{L_i}(a,b)$, contradicting the fact that $(a,b)$ is the closest lattice point to $L_i$ or $\gcd(a,b)=1$. Thus $(a_{\ell},b_{\ell})$ must lie in $\conv\big\{\big(\frac{1}{q},0\big),\big(p - \frac{p-a}{q-b}q,0\big),(a,b)\big\}$ minus the intersection between its boundary with $L_i'$ and $L_i''$. Denote this region $\Delta$.
	
	We parametrize $L_i:qx-py=0$ and $L_i':qx-py-1=0$. The upper bound for the distance of any point $(x,y)\in\Delta$ to $L_i'$ is
	\begin{align*}
		\sup\big\{d_{L_i'}(x,y):(x,y)\in\Delta\big\} &= d_{L_i'}\bigg(p-\frac{p-a}{q-b}q,0\bigg)\\
		&= \frac{1}{\sqrt{p^2+q^2}}\bigg|\bigg(p-\frac{p-a}{q-b}q\bigg)q-1\bigg|\\
		&= \frac{1}{\sqrt{p^2+q^2}}\bigg|\frac{pq(q-b)-(p-a)q^2-q+b}{q-b}\bigg|\\
		&= \frac{1}{\sqrt{p^2+q^2}}\bigg|\frac{q(pq-pb-pq+qa)-q+b}{q-b}\bigg|\\
		&= \frac{1}{\sqrt{p^2+q^2}}\bigg|\frac{b}{q-b}\bigg|.
	\end{align*}
	Consider the cases $b\le \frac{q}{2}$ and $b>\frac{q}{2}$. If the former case holds, then $\big|\frac{b}{q-b}\big|\le1$, which implies $d_{L_i'}(a_{\ell},b_{\ell})\le\frac{1}{\sqrt{p^2+q^2}}$. Let $\delta = d_{L_i'}(a_{\ell},b_{\ell})>0$. Since $b\le\frac{q}{2}$, $2b\le q$. So if we reflect $\Delta$ about the point $(a,b)$, we get the region $\Delta'$, where $\partial\Delta$ and $\partial\Delta'$ are congruent triangles.
	
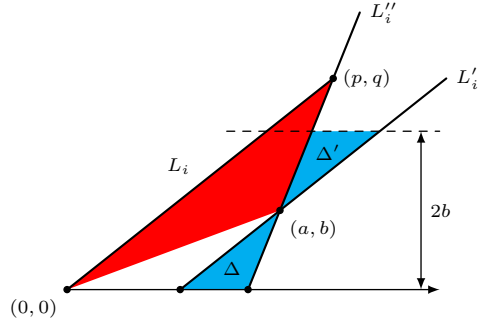
\begin{center}
		\begin{tikzpicture}[mystyle,x=0.7cm,y=0.7cm]
			\fill[red,opacity=0.25] (0,0) -- (5,4) -- (4,1.5) -- cycle;
			\fill[cyan,opacity=0.25] (2.125,0) -- (4,1.5) -- (3.4,0) -- cycle;
			\fill[cyan,opacity=0.25] (8-2.125,3) -- (4,1.5) -- (8-3.4,3) -- cycle;
			
			\draw[->] (0,0) -- (7,0);
			\draw[thick] (0,0) -- node[above left] {$L_i$} (5,4);
			\draw[thick] (2.125,0) -- (7.125,4) node[right] {$L_i'$};
			\draw[thick] (3.4,0) -- (5.5,5.25) node[right] {$L_i''$};
			
			\draw[dashed] (3,3) -- (7,3); 
			\draw[<->] (6.65,0) -- node[right] {$2b$} (6.65,3);
			
			\draw[fill] (0,0) circle (0.04cm) node[below left] {$(0,0)$};
			\draw[fill] (5,4) circle (0.04cm) node[right] {$(p,q)$};
			\draw[fill] (4,1.5) circle (0.04cm) node[below right] {$(a,b)$};
			\draw[fill] (2.125,0) circle (0.04cm);
			\draw[fill] (3.4,0) circle (0.04cm);
			\node at (3.1,0.4) {$\Delta$};
			\node at (8-3.1,3-0.4) {$\Delta'$};
		\end{tikzpicture}
	\end{center}
	Since $(a_{\ell},b_{\ell})\in\Delta$ is a lattice point, $(2a-a_{\ell},2b-b_{\ell})\in\Delta'$ is also a lattice point. But then the distance
	\begin{align*}
		d_{L_i}(2a-a_{\ell},2b-b_{\ell}) = \frac{1}{\sqrt{p^2+q^2}} - \delta < \frac{1}{\sqrt{p^2+q^2}} = d_{L_i}(a,b).
	\end{align*}
	So we have found a lattice point closer to $L_i$ than $(a,b)$, a contradiction. Therefore $(a_{\ell},b_{\ell})\notin\Delta$, which implies it is not possible to append a primitive triangle to $\ell$. 
	
	If $b>\frac{q}{2}$, then $(2a-p,2b-q)$ is a lattice point on the boundary of $\Delta$ that intersects $L_{i}''$. But
	\begin{align*}
		d_{\ell}(2a-p,2b-q) = \frac{|b(2a-p)-a(2b-q)|}{\sqrt{a^2+b^2}} = \frac{|aq-bp|}{\sqrt{a^2+b^2}} = \frac{1}{\sqrt{a^2+b^2}}
	\end{align*}
	is minimized. By uniqueness of the smallest B\'{e}zout coefficients, $(a_{\ell},b_{\ell}) = (2a-p,2b-q)$. Therefore we also cannot append a primitive triangle to $\ell$ in this case because $(a_l, b_l)$ would be collinear with $(a,b)$.
	
	It remains to show it is not possible to append a primitive triangle to the edge $u$. Define $L_i'$ similar to before, but define $L_i''$ to be the line connecting $(0,0)$ and $(a,b)$. By a similar reasoning using convexity to the previous paragraph, we can see that $(a_u,b_u)$ must be contained in $\conv\big\{(a,b),\big(\frac{aq}{b},q\big),\big(\frac{1+pq}{q},q\big)\big\}$ minus its intersection with $L_i'$ and $L_i''$. Again, we denote this region $\Delta$.
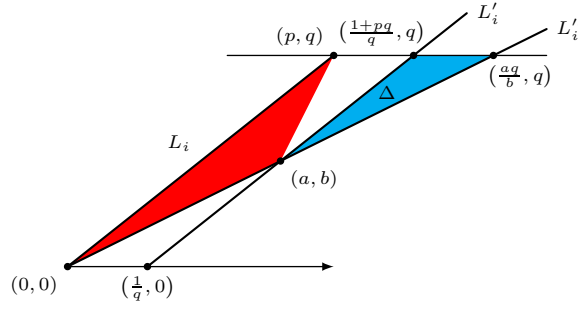
	\begin{center}
		\begin{tikzpicture}[mystyle,x=0.7cm,y=0.7cm]
			\fill[red,opacity=0.25] (0,0) -- (5,4) -- (4,2) -- cycle;
			\fill[cyan,opacity=0.25] (6.5,4) -- (4,2) -- (8,4) -- cycle;
			
			\draw[->] (0,0) -- (5,0);
			\draw (3,4) -- (9,4);
			\draw[thick] (0,0) -- node[above left] {$L_i$} (5,4);
			\draw[thick] (1.5,0) -- (7.5,4.8) node[right] {$L_i'$};
			\draw[thick] (0,0) -- (9,4.5) node[right] {$L_i''$};
			
			\draw[fill] (0,0) circle (0.04cm) node[below left] {$(0,0)$};
			\draw[fill] (5,4) circle (0.04cm) node[above left] {$(p,q)$};
			\draw[fill] (4,2) circle (0.04cm) node[below right] {$(a,b)$};
			\draw[fill] (1.5,0) circle (0.04cm) node[below] {$\big(\frac{1}{q},0\big)$};
			\draw[fill] (6.5,4) circle (0.04cm) node[above,shift={(-0.6,0)}] {$\big(\frac{1+pq}{q},q\big)$};
			\draw[fill] (8,4) circle (0.04cm) node[below,shift={(0.5,0)}] {$\big(\frac{aq}{b},q\big)$};
			\node at (6, 3.3) {$\Delta$};
		\end{tikzpicture}
	\end{center}
	
	The least upper bound for $d_{L_i'}(a_u,b_u)$ is
	\begin{align*}
		\sup\big\{d_{L_i'}(x,y):(x,y)\in\Delta\big\} &= d_{L_i'}\bigg(\frac{aq}{b}, q\bigg)\\
		&= \frac{1}{\sqrt{p^2+q^2}}\bigg|\bigg(\frac{aq}{b}\bigg)q-pq-1\bigg|\\
		&= \frac{1}{\sqrt{p^2+q^2}}\bigg|\frac{q-b}{b}\bigg|.
	\end{align*}
	Consider the cases $b\ge\frac{q}{2}$ and $b<\frac{q}{2}$. If $b\ge \frac{q}{2}$, then $\big|\frac{q-b}{b}\big|\le1$ and $d_{L_i'}(a_u,b_u)\le\frac{1}{\sqrt{p^2+q^2}}$. Again, if we reflect $\Delta$ about $(a,b)$, since $q-2b\ge0$, the reflected $\Delta'\subseteq\conv\big\{(0,0),\big(\frac{1}{n},0\big),(a,b)\big\}$. Therefore $(a_u,b_u)\in\Delta$ is a lattice point, then there must be a lattice point $(2a-a_u,2b-b_u)\in\Delta'$. If we let $\delta = d_{L_i'}(a_u,b_u) > 0$, then $d_{L_i}(2a-a_u,2b-b_u) = d_{L_i}(a,b) - \delta < d_{L_i}(a,b)$, a contradiction. Hence it is not possible to append a primitive triangle to $u$. Finally, if $b<\frac{q}{2}$, then $2b<q$ and $(2a,2b)$ is a lattice point on the boundary $L_i''$ of $\Delta$. We compute
	\begin{align*}
		d_u(2a,2b) = \frac{|2a(q-b)-2b(p-a)-1|}{\sqrt{(p-a)^2+(q-b)^2}} = \frac{1}{\sqrt{(p-a)^2+(q-b)^2}}.
	\end{align*}
	This means that $(2a,2b)$ is the closest point to $u$, and by the uniqueness of the smallest B\'{e}zout coefficients, $(a_u,b_u)=(2a, 2b)$. In this case, it is not possible to append a primitive triangle to $u$.
\end{proof}

This theorem gives an upper bound for the number of times we can add primitive triangles to a polygon before breaking convexity. Although there is no algorithmic way to know the exact number of times we can do so. As soon as we append a triangle to $P$, we change its structure. In other words, an edge by itself may satisfy the conditions to append a primitive triangle, but it may not if either of its adjacent edges had a primitive triangle appended prior.

\begin{example}
	The following are two possible ways to append primitive triangles to a lattice polygon. Note that the number of triangles are different in both cases.
	
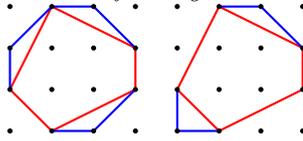
\begin{center}
		\begin{tikzpicture}[mylatstyle]
			\begin{scope}[shift={(1,0)}]
				\draw[thick,red] (0,0) -- (2,1) -- (2,2) -- (0,3) -- (-1,1) -- cycle;
				\draw[thick,blue] (0,0) -- (1,0) -- (2,1);
				\draw[thick,blue] (2,2) -- (1,3) -- (0,3) -- (-1,2) -- (-1,1);
			\end{scope}
			\begin{scope}[shift={(5,0)}]
				\draw[thick,red] (0,0) -- (2,1) -- (2,2) -- (0,3) -- (-1,1) -- cycle;
				\draw[thick,blue] (-1,1) -- (-1,0) -- (0,0);
				\draw[thick,blue] (2,2) -- (1,3) -- (0,3);
			\end{scope}
			\lattice{7}{3}
		\end{tikzpicture}
	\end{center}
\end{example}

\section{Existence and Uniqueness of Convex $n$-gons with $k$ Collinear and Non-Collinear Interior Points}

\subsection{Collinear Points and Half Planes}\label{sec:4.2}

In this final section we show that there does not exist the analogue of $2$-collinear integers for lattice $n$-gons with $n\ge4$. When $n$ is large, it induces the existence of interior points in an $n$-gon with $B=n$. In several cases for $n$ the exact lower bound of the number of interior points is known. For example, in \cite{rabinowitz1,rabinowitz3} Rabinowitz showed a convex lattice $5$-gon must have at least one interior lattice point, and a convex $7$-gon must contain at least four interior points, respectively. More lower bounds can be found in \cite{coleman,rabinowitz2,simpson}. In 1978 Coleman \cite{coleman} conjectured an inequality that relates the number of vertices, boundary, and interior points of any convex lattice $n$-gons. The conjecture is then proved by Ko\l{}odziejczyk and Olszewska (cf. \cite{kolodziejczyk2}).

\begin{theorem}[{\cite[Theorem 15]{kolodziejczyk2}}]
	If $P$ is a convex lattice $n$-gon, then $B(P)\le 2I(P) - n  + 10$.
\end{theorem}

\noindent Given our assumption that $B(P) = n$, we solve for $I$ to get $I(P) \ge n - 5$. Thus for any $n$, we choose $I=k$ large enough so that the bound is satisfied and no contradiction arises.

It suffices to consider collinear points on the $x$-axis because any $k$ collinear points in $\Z^2$ must create a line with rational slope $\frac{p}{q}$. If $p=0$ then we are done, similarly when the points lie on a vertical line then we simply rotate by $\frac{\pi}{2}$. Otherwise we can apply an integral affine transformation to ensure the first point has coordinates $(p,q)$ with $p,q>0$. We then compute the B\'{e}zout coefficients similar to that of Lemmas \ref{lem:existence_convex_m1} and \ref{lem:existence_convex_mn}. Finally, apply another affine map with matrix
\begin{align*}
	A=\begin{pmatrix}
		-b & a\\
		-q & p
	\end{pmatrix}
\end{align*} 
will send the $k$ points to $(1,0),\ldots,(k,0)$.

In \cite{kolodziejczyk1,kolodziejczyk2} the authors defined the \textit{outer hull} of a lattice polygon $P$ to be the closed convex region bounded by lattice lines parallel to the edges of $P$, exterior to $P$, and closest to $P$. They also remarked that the outer hull is well-defined if and only if $P$ is a non-degenerate lattice polygon. Since we cannot define an outer ``bound'' for the $k$ collinear interior points, we define a weaker condition for the possible locations for the vertices of an $n$-gon containing them. Let $R$ denote the union of two open half-planes $\mathcal{H}_1: y > 1$ and $\mathcal{H}_2: y < -1$.

\begin{lemma}\label{lem:upper_bound}
	Let $P$ is a convex $n$-gon with $k\ge3$ collinear interior points. Affine transform $P$ such that the interior points are on the $x$-axis. Then no vertex of $P$ can be an element of $R\cap\Z^2$.
\end{lemma}

\begin{proof}
	Suppose to the contrary that $v\in\partial P$ and $v\in R\cap\Z^2$. Then $v = (x,y)$ with $x\in\Z$ and $y\ge2$ or $y\le-2$. We only prove the case $y\ge2$ since the other case is similar. By convexity, the line between $v$ and $(j,0)$ must be contained in the interior of $P$ for all $j=1,\ldots,k$.
	
	Assume $y=2$. Then, since $k\ge3$, at least one of $x-1,\ldots,x-k$ must be even by the pigeonhole principle. Thus there exists $j\in\{1,\ldots,k\}$ such that $\gcd(x-j,y)=2$. Therefore $L(v,(j,0))>0$, implying the existence of an extra interior point that is not $(1,0),\ldots,(k,0)$, a contradiction.
	
	Assume $y>2$. If there is a $j\in\{1,\ldots,k\}$ such that $\gcd(x-j,y)>1$ then there exists an extra interior point. Thus we assume $y$ is chosen such that $\gcd(x-j,y)=1$ for all $j\in\{1,\ldots,k\}$. Consider the triangle $T=\conv\{(1,0),(k,0),(x,y)\}$. We compute $$\operatorname{area}(T)=\frac{(k-1)y}{2}>k-1.$$ But $B(T) = k + 1$, so by Pick's theorem,
	\begin{align*}
		I(T) = k - 1 - \frac{k + 1}{2} + 1 = \frac{k}{2} - \frac{1}{2}\ge1.
	\end{align*}
	Therefore there exists an extra interior point of $P$, a contradiction.
\end{proof}

\subsection{Results and Conclusions}

\begin{theorem}
	The only integers $n$ in which there is a convex $n$-gon with $n$ boundary and $k\ge3$ interior points in which the interior points are collinear are $3$, $4$, $5$, and $6$.
\end{theorem}

\begin{proof}
	First we show existence of such polygons. For any $k$ collinear points in $\Z^2$ we can affine transform them into $(1,0),\ldots,(k,0)$. Then the $3$-gon $\conv\{(0,0),(k+1,-1),(k,1)\}$ contains exactly $k$ interior points. By Theorem \ref{thm:at_most_n_triangles}, we can append at most $3$ primitive triangles to obtain a $4$-, $5$-, and $6$-gon containing the same number of interior points. Below is an example of appending three primitive triangles for $k=5$.
	
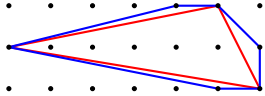
\begin{center}
		\begin{tikzpicture}[mylatstyle]
			\draw[thick,red,shift={(0,1)}] (0,0) -- (6,-1) -- (5,1) -- cycle;
			\draw[thick,blue,shift={(0,1)}] (0,0) -- (5,-1) -- (6,-1) -- (6,0) -- (5,1) -- (4,1) -- cycle;
			\lattice{6}{2}
		\end{tikzpicture}
	\end{center}
	
	 Now we show $6$ is the best we could do. Assume for contradiction that we can construct an $n$-gon, $n\ge7$ with $k$ collinear interior points. By Lemma \ref{lem:upper_bound}, the vertices of said $n$-gon must be contained in $\conv\{(x_1,-1),(x_2,-1),(x_2,1),(x_1,1)\}$ with $x_1\le0$ and $x_2\ge k+1$. Note that we can only have at most $2$ vertices on the lines $y=\pm1$ because otherwise we would have collinear boundary points. Then there are $n-4 \ge 3$ remaining vertices. Then there must be at least two lattice points with $x$-component both less than $0$ or greater than $k+1$. Further, their $y$-component cannot be the same because it will break convexity. Hence at least one of the two points must have $y$-component $\pm1$, which induces at least three collinear boundary points.
\end{proof}

\begin{theorem}
	For any $k\ge3$ there exists a convex $n$-gon, $n\in\{4,5,6\}$ such that the $k$ interior points are not collinear.
\end{theorem}

\begin{proof}
	For each $k\ge3$ define the $4$-gon
	\begin{align*}
		P_k = \begin{cases}
			\conv\big\{(0,0),\big(\frac{k}{2},-1\big),\big(\frac{k}{2}+1,1\big),(1,2)\big\} &\text{if }k\text{ even},\\
			\conv\big\{(0,0),\big(\frac{k+1}{2}+1,-1\big),\big(\frac{k+1}{2},1\big),(1,2)\big\} &\text{if }k\text{ odd}.
		\end{cases}
	\end{align*}
	For example, below are $P_6$ and $P_7$.
	
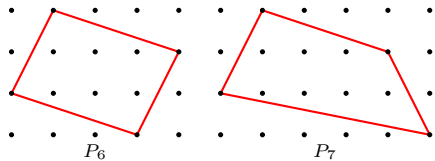
\begin{center}
		\begin{tikzpicture}[mylatstyle]
			\draw[thick,red,shift={(0,1)}] (0,0) -- (3,-1) -- (4,1) -- (1,2) -- cycle;
			\draw[thick,red,shift={(5,1)}] (0,0) -- (5,-1) -- (4,1) -- (1,2) -- cycle; 
			\node[below] at (2,0) {$P_6$};
			\node[below] at (7.5,0) {$P_7$};
			\lattice{10}{3}
		\end{tikzpicture}
	\end{center}
	Then $\{P_k\}$ is a sequence of $4$-gons with $k$ non-collinear interior points. By Theorem \ref{thm:at_most_n_triangles}, we can append at most four primitive triangles to each $P_k$, but we only need to append to two edges. Below are examples for $P_6$ and $P_7$.
	
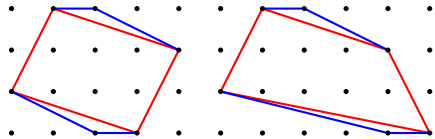
\begin{center}
		\begin{tikzpicture}[mylatstyle]
			\draw[thick,red,shift={(0,1)}] (0,0) -- (3,-1) -- (4,1) -- (1,2) -- cycle;
			\draw[thick,red,shift={(5,1)}] (0,0) -- (5,-1) -- (4,1) -- (1,2) -- cycle;
			\draw[thick,blue] (0,1) -- (2,0) -- (3,0);
			\draw[thick,blue] (1,3) -- (2,3) -- (4,2);
			\draw[thick,blue] (5,1) -- (9,0) -- (10,0);
			\draw[thick,blue] (6,3) -- (7,3) -- (9,2);
			\lattice{10}{3}
		\end{tikzpicture}
	\end{center}
\end{proof}

In fact, if $k\notin\{3,5\}$ we can apply Lemmas \ref{lem:existence_convex_m1} and \ref{lem:existence_convex_mn} to the sequence $\{P_k\}$ from above and go all the way up to $8$-gons containing $k$ non-collinear interior points. But since $n$-gons having $k$ collinear interior points only go up to $n=6$, we only need to go up to $6$ to show that there does not exist $2$-collinear integers for $n$-gons, $n\ge4$. The reason why $k=3$ does not work is due to \cite{rabinowitz3} and $k=5$ does not work is due to \cite[Theorem 4.2]{kolodziejczyk1}.

\section*{Acknowledgements} 

This research was generously supported by the William and
Linda Frost Fund in the Cal Poly Bailey College of Science and
Mathematics. This material is based upon work supported by the National Science Foundation under Grant No. 2015553. The authors thank the reviewers for their time and efforts in reviewing the paper. 


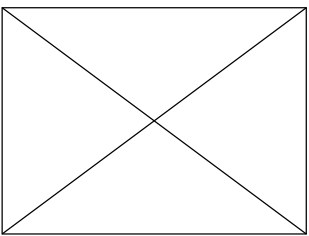
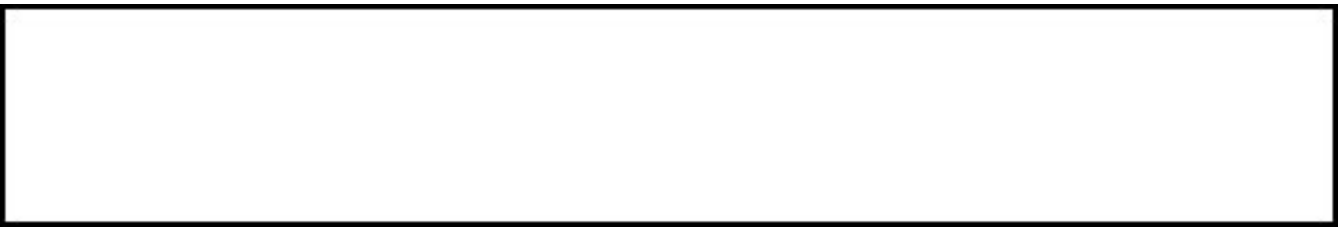
\end{document}